\documentclass[11pt,reqno]{amsart}

\usepackage[utf8]{inputenc}

\usepackage{a4wide}
\usepackage[italian,english]{babel}

\usepackage{fancyhdr}
\usepackage{indentfirst}
\usepackage{graphicx}
\usepackage{newlfont} 
\usepackage{dsfont}
\usepackage{amssymb}
\usepackage{amsmath}
\usepackage{latexsym}
\usepackage{amsthm}

\usepackage{amsthm}

\usepackage{hyperref}
\usepackage{faktor}
\usepackage{mathtools}
\usepackage{thmtools}
\usepackage{tikz}
\usetikzlibrary{matrix,calc}
\usetikzlibrary{cd}
\usepackage{cleveref}

\newtheorem{thm}{Theorem}[section]
\newtheorem{defin}[thm]{Definition}
 
\newtheorem{cor}[thm]{Corollary}
\newtheorem{lemma}[thm]{Lemma}
\newtheorem{prop}[thm]{Proposition}

\theoremstyle{remark}
\newtheorem*{rmk}{Remark}

\newcommand{\Q}{\mathbb{Q}}

\newcommand{\Z}{\mathbb{Z}}

\newcommand{\op}{\operatorname}

\newcommand{\mc}[1]{\mathcal{[1]}}

\newcommand{\R}{\mathbb{R}}

\newcommand{\img}{\operatorname{Im}}

\let\ForAll\forall
\renewcommand{\forall}{\quad \ForAll}

\newcommand{\dual}[1]{#1\,\check{\vrule height1.3ex width0pt}}

\newcommand{\LambdaY}{\Lambda_{Nik}}
\newcommand{\LambdaX}{\Lambda_{K3^{[2]}}}

\newcommand{\alert}[1]{\par \textcolor{red}{\textbf{\huge{#1}}}\par}
\renewcommand{\alert}[1]{}

\usepackage{biblatex}
\bibliography{bibliography}

\author{Giacomo Nanni}
\address{Dipartimento di Matematica \\
	Universit\`a di Bologna\\
	Piazza di Porta San Donato 5\\
	40127 Bologna, Italy}
\email[G.~Nanni]{giacomo.nanni13@unibo.it}

\makeatletter

\@namedef{subjclassname@2020}{\textup{2020} Mathematics Subject Classification}
\subjclass[2020]{Primary 14J42; Secondary 14D05}

\keywords{Symplectic orbifolds, Monodromy, Nikulin orbifolds}
\thanks{The author is part of INdAM
research group “GNSAGA” and was partially supported by it. The author acknowledges the support of the European Union - NextGenerationEU under the National Recovery and Resilience Plan (PNRR) - Mission 4 Education and research - Component 2 From research to business -Investment 1.1 Notice Prin 2022 - DD N. 104 del 2/2/2022, from title "Symplectic varieties: their interplay with Fano manifolds and derived categories", proposal code 2022PEKYBJ – CUP J53D23003840006. The author was also supported by the DFG through the research grant Le 3093/5-1.}

\title{Monodromy of Nikulin orbifolds}

\begin{document}

\maketitle
\vspace{-25pt}
\begin{abstract}
    We give a new proof for the maximality of the monodromy group of a Nikulin orbifold, a symplectic orbifold arising as terminalisation of a symplectic quotient of a \text{$K3^{[2]}$-type} fourfold.
\end{abstract}
\vspace{-5pt}
\section{Introduction}
A Nikulin orbifold is an irreducible symplectic orbifold (see \cite[Definition 3.1Proposition 2.8]{MenetGlobalTor}) obtained as the $\Q$-factorial terminalisation of the quotient of a $K3^{[2]}$-type fourfold by a symplectic involution. More concretely, let $X$ be a $K3^{[2]}$-type fourfold, let $\iota\in \op{Aut}(X)$ a symplectic involution. The quotient $\hat X:=\faktor{X}{\iota}$ is singular along the disjoint union of a K3 surface $\Sigma$ and 28 isolated singular points. The $\Q$-factorial terminalisation is obtained as the blow-up $Y:=Bl_\Sigma \hat X\rightarrow \hat X$ of the K3 surface $\Sigma$.

The second cohomology group of a Nikulin orbifold is endowed with a Beauville-Bogomolov-Fujiki (BBF) form \cite[Theorem 3.17]{MenetGlobalTor} for which the lattice structure has been computed in \cite{MenetBBFNik}
 to be isomorphic to:
 \begin{equation}\label{eq:LambdaY}
   \Lambda_{N}:=U(2)^{ 3}\oplus E_8(-1)\oplus\langle -2 \rangle^{ 2}  
 \end{equation}
 where $U$ is the hyperbolic lattice, $E_8$ is the rank eight even unimodular positive definite lattice and $\langle -2 \rangle$ is a rank one lattice with primitive generator of square $-2$.

For $Y$ a Nikulin orbifold, let $Mon^2(Y)$ be the monodromy group of $Y$ and let $O^+(Y)$ be the
group of orientation-preserving isometries (with respect to the real spinor-norm). The monodromy group has been recently computed in \cite[Theorem 1.1]{StevellMonNik}.

\begin{thm}\label{fullMon}
    Let $Y$ be a Nikulin orbifold. Then $Mon^2(Y)=O^+(Y)$.
\end{thm}
We give an alternative proof which makes the underlying geometry more apparent.
While the approach by Brandhorst, Menet and Muller is based on producing a set of reflections of positive spinor norm and then proving by lattice-theoretic means that it generates the entire monodromy, ours derives the result from the analogous one on the corresponding $K3^{[2]}$ fourfold relying on the monodromy operators introduced in \cite{monNik}. Moreover, it has the advantage of relying on more elementary lattice theory.

\section{Preliminaries}

\subsection{Lattice notation}\label{latticeStuff}
 In what follows, for a lattice $\Lambda$ and $v,w\in \Lambda$ we use the notation $(v,w)_\Lambda$ to denote their product in $\Lambda$. When there is no confusion possible, we may omit the lattice from the notation by writing $(v,w)$. For a direct sum of lattices $\Lambda\oplus\Gamma$ we will write $\lambda\oplus\gamma$ to denote the element with components $\lambda\in\Lambda,\gamma\in\Gamma$. We will denote by $\dual \Lambda:=\op{Hom}_\Z(\Lambda,\Z)$ the dual of $\Lambda$ and $A_\Lambda:=\faktor{\dual{\Lambda}}{\Lambda}$ the discriminant group.  
 
 The group of isometries of $\Lambda$ will be denoted $O(\Lambda)$. We denote by $R_v$ the reflection around a non-isotropic vector $v\in \Lambda$ defined as  $R_v(\alpha)=\alpha-2 \frac{(\alpha,v)}{(v,v)} v$. The real spinor norm is the character $O(\Lambda)\rightarrow \{\pm1\}$ such that $R_v\mapsto -sign(v^2)$. We say that an isometry is orientation-preserving if it lies in the kernel $O^+(\Lambda)$ of the real spinor norm. 
 \begin{rmk}
     Notice that the definition of orientation-preserving in this context does not coincide with the usual one in terms of the sign of the discriminant. This possibly unfortunate but established terminology is motivated in \cite[Section 4]{MarkmanSurvey} and in agreement with \cite{StevellMonNik}. 
 \end{rmk}

The identification $\dual\Lambda\cong\{v\in \Lambda\otimes\Q, (v,\Lambda)\subset\Z\}$ induces a $\Q$-valued quadratic form defined on $\dual\Lambda$ such that the embedding $\Lambda\hookrightarrow\dual \Lambda$ is an isometry. Moreover, if $\Lambda$ is an even lattice, this induces a quadratic form on the discriminant group $A_\Lambda$ taking values in $\faktor{\Q}{2\Z}$. We will denote $O(A_\Lambda)$ the subgroup of group homomorphisms preserving the quadratic form. An isometry $\psi\in O(\Lambda)$ induces an isometry $\bar\psi\in O(A_\Lambda)$ by $\bar\psi(\alpha+\Lambda)=\alpha\circ \psi^{-1}+\Lambda$ for $\alpha\in \dual\Lambda$.

\subsection{Monodromy}
Let $(Y,\eta)$ be a marked Nikulin orbifold, which is to say: a Nikulin orbifold $Y$ together with an isomorphism of lattices $\eta:(H^2(Y,\Z),q_Y)\rightarrow\LambdaY$. 
Let $Mon^2(Y)\leq O(H^2(Y,\Z))$ be the subgroup of monodromy operators in the group of isometries of the BBF lattice of $Y$. It is a general fact %{\bf{(I may cite \cite{StevellMonNik}, as they state it in the introduction, but probably I should find a better refernce)}},
that monodromy operators are orientation-preserving (see for example \cite[Section 4]{MarkmanSurvey}, or the introduction of \cite{StevellMonNik}). Therefore $Mon^2(Y)\leq O^+(H^2(Y,\Z))$.
%The main result of this section is

 Let $O(\LambdaY)$ be the group of isometries of the lattice $\LambdaY$ and $O^+(\LambdaY)\leq O(\LambdaY)$ the subgroup of orientation preserving isometries. For a marked Nikulin orbifold, we denote $Mon(\LambdaY,\eta):=\eta\circ Mon(Y,\Z)\circ\eta^{-1}\leq O^+(\LambdaY)$. This, a priori, may depend on $\eta$. 
 
 We fix $\eta$ as in \cite[Section 2]{MenetBBFNik}. Therefore, we will drop $\eta$ from the notation. We briefly recall its description. Let $X$ a $K3^{[2]}$ fourfold, $\iota$ a symplectic involution on $X$ and $\Sigma$ such that $Y=Bl_\Sigma(\faktor{X}{\iota})$. Push-forward to the quotient and pull-back through the blow-up gives an isometric embedding of a twist of the $\iota$-invariant part of the BBF lattice  $H^2(X,\Z)^{\iota^*}(2)\rightarrow H^2(Y,\Z)$. Moreover, denoting $E_\Sigma\in H^2(Y,\Z)$ the exceptional divisor over $\Sigma$, $H^2(X,\Z)^{\iota^*}(2)\oplus \Z E_\Sigma$ is a finite index sublattice of $H^2(Y,\Z)$.  Denote \begin{equation}\label{eq:LambdaX}
     \LambdaX:= U^{3}\oplus E_8(-1)^{2}\oplus \langle -2 \rangle.
 \end{equation} Denote $\delta_X$ the generator of the $\langle -2 \rangle$ summand. A marking $\eta_X$ of $X$ can always be chosen in such a way that the involution induced by $\iota^*$ on $\LambdaX$ acts exchanging the two $E_8$ summands (as in \cite[Proposition 2.6]{MenetBBFNik}, based on \cite[Theorem 1.3]{MongardiInvolutionsK3[2]}). The fixed locus of this action is isomorphic to \begin{equation}\label{eq:LambdaFix}
     \Lambda_{fix}:=U^3\oplus E_8(-2)\oplus \langle -2 \rangle
 \end{equation} which embeds in $\LambdaX$ as $u\oplus e\oplus k\delta_X\mapsto u\oplus e\oplus e\oplus k\delta_X$ (where $u\in U^3, e\in E_8,k\in \Z$).
 Moreover, $\Lambda_{fix}(2)$ embeds in $\LambdaY$. To see this explicitly, define $\delta_Y,\Sigma_Y\in \LambdaY$ such that $\frac{\delta_Y+\Sigma_Y}{2},\frac{\delta_Y-\Sigma_Y}{2}$ are the generators of the two $\langle -2 \rangle$ summands. Then, $u\oplus e\oplus k\delta_X\mapsto (u\oplus 2e\oplus k\frac{\delta_Y+\Sigma_Y}{2}\oplus k\frac{\delta_Y-\Sigma_Y}{2})$ is an isometric embedding $\Lambda_{fix}(2)\rightarrow\LambdaY$. The main result of \cite{MenetBBFNik} proves that the isometry $H^2(X,\Z)^{\iota^*}(2)\oplus \Z E_\Sigma\rightarrow \Lambda_{fix}(2)\oplus \Z\Sigma_Y$ extends to a marking $H^2(Y,\Z)\rightarrow\LambdaY$. 
 
\subsection{Sum of twisted unimodular lattices}
For a lattice $\Lambda$ let $Rays(\Lambda):=\{\R_{\geq 0} v\subset \Lambda\otimes \R,v\in \Lambda\}$ be the set of rays.
We give a slightly more explicit statement of \cite[Lemma 6.2]{monNik}
\begin{prop}
    Let $N,M$ unimodular lattices. Let $\phi:N\oplus M(2)\hookrightarrow N(2)\oplus M$ the map $n\oplus m\mapsto n\oplus 2m$. Then: 
    \begin{itemize}
        \item The map $\Phi:O(N(2)\oplus M)\rightarrow O(N\oplus M(2))$ defined as $f\mapsto \phi^{-1} \circ f\circ\phi$ is well defined and is an isomorphism of groups.
        \item The map $\rho: Rays(N\oplus M(2))\rightarrow Rays(N(2)\oplus M)$ defined as $\R_{\geq 0} v\mapsto \R_{\geq 0}\phi(v)$ is bijective.
    \end{itemize}
    \proof
    First, observe that the image $\img \phi$ is the sublattice of elements of even divisibility. Therefore it is preserved by any isometry, proving the well-posedness $\Phi:O(N(2)\oplus M)\rightarrow End(N\oplus M(2))$. To prove $\img \Phi\subset O(N\oplus M(2))$, one can observe that \begin{align*}
        (\phi(n\oplus m),\phi(n'\oplus m'))_{N(2)\oplus M}=(n\oplus 2m,n'\oplus 2m')_{N(2)\oplus M}=(n,n')_{N(2)}+(2m,2m')_{M}=\\2(n,n')_N\oplus4(m,m')_M=2(n\oplus m,n'\oplus m')_{N\oplus M(2)}.
    \end{align*} This concludes, as $(\Phi(f)(v),\Phi(f)(w))=\frac{1}{2}(f(\phi(v)),f(\phi(w)))=(v,w)$.

    For the second part, this automatically follows from the fact that tensoring with $\R$, $\phi$ extends to an isomorphism of $\R$-vector spaces. \endproof
\end{prop}
\begin{rmk}
    From the definitions, it immediately follows that $\phi(\Phi(f)(v))=f(\phi(v))$.
\end{rmk}
Following Menet-Riess, we introduce the lattice \begin{equation}\label{eq:Lambda1}
    \Lambda_1:=U^3\oplus E_8(-2)\oplus \langle -1 \rangle^2.
\end{equation} We denote $\delta_1,\Sigma_1$ the two -2 classes such that $\frac{\delta_1+\Sigma_1}{2},\frac{\delta_1-\Sigma_1}{2}$ are the generators of the $\langle -1 \rangle$ summands. The previous proposition can be applied by choosing $N=U^3\oplus \langle -1 \rangle^2$ and $M=E_8(-1)$. We get $\phi:\Lambda_1\rightarrow \LambdaY$ and all the corresponding maps between rays and isometries.
\subsection{Isometries of the invariant lattice}
Let $\sigma\in O(\LambdaX)$ denote the involution exchanging the two $E_8$ factors. 
%As stated before, the invariant lattice $\Lambda$
%Let $\Lambda$ be a lattice and $\sigma\in O(\Lambda)$ an isometry. 
Let $\Lambda^\sigma,\Lambda_{\sigma}$ denote the invariant and coinvariant lattice of $\sigma$. Let $O(\LambdaX)^\sigma$ be the group of $\sigma$-equivariant isometries of $\LambdaX$. Then restriction to $\Lambda^\sigma$ gives a group morphism $O(\LambdaX)^\sigma\rightarrow O(\Lambda^\sigma)$. We prove that this map is surjective.

%Since the direct sum $\Lambda^\sigma\oplus\Lambda_\sigma\subset \Lambda$ is a full rank sublattice lattice,
%Let $\psi\in O(\Lambda_\fix)$

We start by recalling the following general criterion: let $L$ be a lattice and $T\subset\Lambda$ a primitive sublattice. By the sequence of embeddings $$T\oplus T^\perp\hookrightarrow L\subset \dual{L}\hookrightarrow \dual{(T\oplus T^\perp)}\cong \dual{(T)}\oplus\dual{(T^\perp)} $$ we define a subgroup $M:=\faktor{L}{(T\oplus T^\perp)}\subset A_{T}\oplus A_{T^\perp}$. Let $p,p_\perp$ be the canonical projections of $A_{T}\oplus A_{T^\perp}$. We define $\overline M$ (resp.$\overline M_\perp$) the image of $M$ through $p$ (resp. $p_\perp$) and $\gamma=p\circ {p_\perp^{-1}}_{|\overline M_\perp}:\overline M_\perp\rightarrow \overline M$. Then the following  holds:
%By the sequence of embeddings $$L^\sigma\oplus\Lambda_\sigma\hookrightarrow \Lambda\subset \dual{\Lambda}\hookrightarrow \dual{(\Lambda^\sigma\oplus\Lambda_\sigma)}\cong \dual{(\Lambda^\sigma)}\oplus\dual{(\Lambda_\sigma)} $$ we define a subgroup $M:=\faktor{L}{(\Lambda^\sigma\oplus\Lambda_\sigma)}\subset A_{\Lambda^\sigma}\oplus A_{\Lambda_\sigma}$. Let $p^\sigma,p_\sigma$ be the canonical projections of $A_{\Lambda^\sigma}\oplus A_{\Lambda_\sigma}$. We define $M^\sigma,M_\sigma$ the image of $M$ through these projections and $\gamma=p^\sigma\circ p_\sigma^{-1}_{|M_\sigma}:M_\sigma\rightarrow M^\sigma$. Then the following criterion holds
\begin{prop}\cite[Corollary 1.5.2]{nikulin}\label{extensionCriterion}
    Let $L$ be an even lattice, $T\subset L$ a primitive sublattice and $\phi\in O(T)$. Then there exists an isometry $\Phi\in O(L)$ such that $\Phi_{|T}=\phi$ if and only if there exists $\psi\in O(T^\perp)$ such that $\gamma\circ\bar\psi=\bar\phi\circ\gamma$.
\end{prop}
To apply the criterion to our case we will need the following computation:
\begin{lemma}\label{isomorphismOfDiscr}
    The discriminant group $A_{E_8(-2)}$ is isomorphic to $\faktor{E_8}{2E_8}$. Moreover, as quadratic modules ${A_{E_8(-2)}\cong\faktor{E_8}{2E_8}(2)}$.%The quadratic form on the discriminant group coincides with the one induced by $E_8$.
    \proof
    Recall that as groups $E_8$ and $E_8(-2)$ are the same and therefore $\dual E_8= \dual {E_8(-2)}$ as groups. %To avoid confusion, in this proof we  denote with $(,)_2$ the product of $E_8(-2)$ and with $(,)$ the one of $E_8$.

    Since $E_8$ is unimodular, the map $E_8\rightarrow \dual E_8=\dual {E_8(-2)}$ sending $e\mapsto (e,-)_{E_8}=-\frac{1}{2}(e,-)_{E_8(-2)}$ is an isomorphism. As the preimage of the sublattice $E_8(-2)\subset \dual {E_8(-2)}$ is $2E_
    8$ we get the group isomorphism $\alpha:\faktor{E_8}{2E_8}\rightarrow A_{E_8(-2)}$.

    Computing the quadratic form of $A_{E_8(-2)}$ for ${e+2E_8\in \faktor{E_8}{2E_8}}$ under this identification gives  \begin{align*}
(\alpha(e+2E_8),\alpha(e+2E_8))_{A_{E_8(-2)}}=(\frac{1}{2}(e,-)_{E_8(-2)},\frac{1}{2}(e,-)_{E_8(-2)})_{\dual {E_8(-2)}\otimes\Q}+2\Z=\\=\frac{1}{4}(e,e)_{E_8(-2)}+2\Z=\frac{1}{2}(e,e)_{E_8}+2\Z\in \faktor{\Z}{2\Z} . 
    \end{align*} 
    \endproof
    
    % Consider the map $E_8\rightarrow A_{E_8(-2)}$ sending $e\rightarrow \frac 1 {\op{div}_{E_8(-2)} e}[(e,-)_{E_8(-2)}]$. As primitive elements of $E_8$ have divisibility $2$, we can rewrite $\frac 1 {\op{div}_{E_8(-2)} e}(e,-)_{E_8(-2)}=\frac{1}{\op{div}_{E_8} e}([e,-)_{E_8}]$. This is surjective and it has kernel $2E_8$.
\end{lemma}%fix $T=\Lambda^\sigma\subset \Lambda$. 
    We briefly recall the notion of characteristic vector, as we will need it in the following proof.
    \vspace{-5pt}
    \begin{defin}\label{def:characteristicVector}
         Let $V$ a vector space over a field of characteristic 2 and let $b:V\rightarrow V$ a non-degenerate bilinear form. Since the map $v\mapsto b(v,v)$ for $v\in V$ is linear in characteristic 2, there exist a unique vector $w\in V$ such that $b(w,v)=b(v,v)$ for all $v\in V$. We call $w$ the characteristic vector.
    \end{defin}

\begin{prop}\label{extendIsometries}
    Every isometry $\phi\in O(\Lambda^\sigma)$ extends to an isometry $\Phi\in O(\LambdaX)$.
%such that $\Phi_{|\Lambda^\sigma}=\phi$.
    \proof
    We keep the notation $T:=\Lambda^\sigma\subset \LambdaX$.
    Then $A_T\cong A_{E_8(-2)}\oplus A_{\langle 2\rangle}$ and $A_{T^\perp}\cong A_{E_8(-2)}
$. Moreover, a direct computation shows that  $M\subset A_T\oplus A_{T^\perp}$ (as defined at the beginning of the section) is the sublattice of elements of the form $[v]_{A_{E_8(-2)}}\oplus0\oplus [-v]_{A_{E_8(-2)}}\in A_{E_8(-2)}\oplus A_{\langle 2\rangle}\oplus A_{E_8(-2)}$. Therefore, the map $\gamma:\overline{M}_\perp\rightarrow\overline{M}$ is the multiplication by $-1$ between two copies of $A_{E_8(-2)}$.
    By \cref{extensionCriterion}, we only have to prove that $\bar \phi$ preserves the $A_{E_8(-2)}$ component of $A_T$ and that there exists an isometry $\psi\in O(T^\perp)$ such that $-\bar\psi=-\bar\phi$.
    
    For the first part, notice that the generator of the $A_{\langle-2\rangle}$ component is the characteristic vector of $A_T$ (as in \cref{def:characteristicVector}). By its unicity, it follows that any isometrie preserves it. Therefore, its orthogonal complement $A_{E_8(-2)}\subset A_T$ is also preserved.
    
    For the second part, it is enough to show that the map ${O(E_8(-2))\rightarrow O(A_{E_8(-2)})}$ is surjective, as $T^\perp\cong E_8(-2)$. 
By \cref{isomorphismOfDiscr}, this is the same as to study the map ${O(E_8)\rightarrow O(\faktor{E_8}{2E_8})}$. But this is surjective by \cite[Exercise VI.4.1]{BourbakiLie456}.  \endproof
\end{prop}
\begin{rmk}
    By the proof of the criterion \cref{extensionCriterion}, using the same notation of the statement, the restriction of the extension $\Phi$ to $T\oplus T^\perp$ can be assumed to be of the form $\phi\oplus\psi$. Since in the discriminat group $\op{id}_{A_{T^\perp}}=-\op{id}_{A_{T^\perp}}$, we can assume $\Phi\in O^+(\LambdaX)$ by replacing $\psi$ with $-\psi$ if needed.
\end{rmk}
\begin{cor}\label{cor:equivariantExtension}
    The restriction map $O(\LambdaX)^\sigma\rightarrow O(\Lambda^\sigma)$ is surjective.
    \proof
    By \cref{extendIsometries}, any $\phi\in O(\Lambda^\sigma)$ extends to $\Phi\in O(\Lambda)$. We only need to prove that $\Phi$ is $\sigma$-equivariant.
    This is trivial on $\Lambda^\sigma\oplus {\Lambda^\sigma}^\perp$ as the action of $\sigma$ on each one of the summands is respectively $\op{id}_{\Lambda^\sigma},-\op{id}_{{\Lambda^\sigma}^\perp}$. As $\Lambda^\sigma\oplus {\Lambda^\sigma}^\perp$ is an index 2 sublattice for any $v\in \Lambda$ we have $2v\in \Lambda^\sigma\oplus {\Lambda^\sigma}^\perp$. Therefore $$\sigma\phi(v)=\frac{1}{2}\sigma\phi(2v)=\phi(\sigma v)$$
    \endproof
\end{cor}

\section{Results}
\subsection{Preservation of $\Sigma_Y$ under isometries}
Here, we review the classes of monodromy operators introduced in \cite{monNik}.
Let $f\in Mon(\LambdaX)$ be a $\iota^*$-equivariant isometry. Then $f_{|\Lambda_{fix}}\in O(\Lambda_{fix})=O(\Lambda_{fix}(2))$. Define $\hat f\in O(\LambdaY)$ as the unique isometry such that $\hat f(\Sigma_Y)=\Sigma_Y$ and $\hat f_{|\Lambda_{fix}(2)}=f_{|\Lambda_{fix}}$. By \cite[Corollary 6.4]{monNik}, $\hat f\in Mon(\LambdaY)$. We say $\hat f$ is a monodromy operator induced by the smooth 4fold.
Also, reflections $R_v$ around (-4) and (-2) classes $v$ of divisibility 2 %, defined by $R_v(\alpha)=\alpha-2 \frac{(\alpha,v)}{(v,v)} v$ 
lie in $Mon(\LambdaY)$ (\cite[Corollary 6.18]{monNik}).

We will denote $G\leq Mon(\LambdaY)$ the subgroup generated by these two types of operators. The orbits of primitive vectors in $\LambdaY$ for the action of $G$ are classified in \cite[Theorem 6.15]{monNik}.%We will define $O_{j,i}\subset \LambdaY$ to be the $i$-th orbit with index $i$ (so that for example in their notation $L_i\in O_{1,i}$) following their enumeration in the statement of \cite[Theorem 6.15]{monNik}.
\begin{rmk}
    Notice that for any vector the condition of spanning a ray of divisibility 1 in $\Lambda_ 1$ is preserved by isometries. In fact, if $v\in\LambdaY$ verifies this condition it means $\R v=\R\phi(w)$ for some $w\in\Lambda_1$ with $\op{div} w=1$. But then taking an isometry $f\in O(\LambdaY)$ one gets $\R f(v)=f(\R v)=f(\R\phi(w))=\R\phi(\Phi(f)w)$.
\end{rmk}
\begin{lemma}
The monodromy orbit of $\Sigma_Y$ is preserved under isometries.
    \proof
    
    Consider $f\in O(\LambdaY)$. Let us denote $O_1,O_2,O_3$ respectively the $G$-orbits of (in the notation of \cite[Theorem 6.15]{monNik}) $L_{-1}^{(2)},-\delta',2e_2^{(1)}-\delta'$.
    Since $f$ is an isometry, $f(\Sigma_Y)$ is primitive with square -4 and divisibility 2.   Therefore, we deduce $f(\Sigma_Y)\in O_{1}\cup O_{2}\cup O_{3}$. We show now that $f(\Sigma_Y)\in O_{2}$. 
    
    If $f(\Sigma_Y)\in O_{1}$, then $f(\Sigma_Y)$ would correspond to a ray of divisibility 1 in $\Lambda_1$. But this is a condition preserved by isometry. Since it was not the case for $\Sigma_Y$ (as the component of $\Sigma_Y$ on the $U(2)^3$ factor is divisible by 2), $f(\Sigma_Y)\not\in O_{3}$.
    Now we prove that the projection $\overline{f(\Sigma_Y)_{E_8}}$ of the $E_8$ component $f(\Sigma_Y)_{E_8}$ to $\faktor{E_8}{4E_8}$ is zero, concluding therefore that $f(\Sigma_Y)\in O_{2,0}$.

    From the explicit definition of $\phi$, we have $\phi(\Sigma_1)=\Sigma_Y$. The projection $\overline{\Sigma_1}$ of $\Sigma_1$ on the discriminant group $\faktor{E_8}{2E_8}$ of $\Lambda_1$ is 0. Denote $\overline{\Phi(f)}$ the isometry of the discriminant group induced by $\Phi(f)$. Then $0=\overline{\Phi(f)}(\overline{\Sigma_1})=\overline{(\Phi(f)(\Sigma_1))}$. This means that the $E_8$ component $(\Phi(f)(\Sigma_1))_{E_8}$ is even. But then $(f(\Sigma_Y))_{E_8}=(\phi(\Phi(f)(\Sigma_1)))_{E_8}=2(\Phi(f)(\Sigma_1))_{E_8}$ is divisible by 4, which is what we wanted. \endproof
\end{lemma}
\subsection{Proof of the main result}
\begin{proof}[Proof of \cref{fullMon}]
%Let $f\in O^+(\LambdaY)$ be an isometry. 
Let $f\in O(\LambdaY)$ be an isometry. We show that either $f$ or $-f$ is a monodromy operator. 

By the previous lemma, there exists an element $g\in G$ such that $g\circ f(\Sigma_Y)=\Sigma_Y$. To prove that $f$ is a monodromy operator is the same as to prove $g\circ f$ is one. Therefore we can assume from the start that $f$ fixes $\Sigma_Y$. 
It follows that $f$ preserves $\Sigma_Y^\perp$. But $\Sigma_Y^\perp$ is an overlattice of $\Lambda_{fix}(2)\subset\LambdaY$ (cf. \cref{eq:LambdaFix}). Moreover, $\Lambda_{fix}(2)$ is the sublattice of elements of even divisibility inside $\Sigma_Y^\perp$. Indeed, the condition for an element  $u\oplus e\oplus k\oplus k'\in \LambdaY$ (where we use the same description as in ) to be orthogonal to $\Sigma^\perp$ gives $0=2k-2k'$, and the parity of its divisibility is the same as the parity of the divisibility of $e$: taken another element $\tilde u\oplus \tilde e\oplus \tilde k\oplus \tilde k\in \Sigma^\perp$ their pairing is $(u,\tilde u)_{U^3(2)}+(e,\tilde e)_{E(-1)}-4k\tilde k\equiv_{\op{mod} 2}(e,\tilde e)_{E(-1)}$. From the unimodularity of $E(-1)$, an element $e\in E(-1)$ has even divisibility if and only if it is divisible by 2.  %square divisible by 4 inside $\Sigma_Y^\perp$. In fact, the condition for an element  $u\oplus e\oplus k\oplus k'\in \LambdaY$ (where we use the same description as in ) to be orthogonal to $\Sigma^\perp$ gives $0=2k-2k'$, and the square is $(u,u)_{U(2)}+(e,e)_{E(-1)}+-2k-2k'$. Since $U$ is even, the squares of elements in $U(2)$ are multiples of 4, and we deduce that $k=k'$ and $0\equiv_{\op{mod} 4}(e,e)_{E(1)}$ which implies $e$ to be divisible by 2. 

It follows that $\Lambda_{fix}(2)\subset\LambdaY$ is preserved by $f$.
By \cref{cor:equivariantExtension} each isometry of $\Lambda_{fix}$ can be extended to a $\iota^*$-equivariant isometry of $\Lambda_{X}$. Therefore we can consider the isometry $\hat f\in O(\LambdaX)$ extending $f_{|\Lambda_{fix}}$.
Since $O^+(\Lambda_{X})=Mon(X)$ then either $\hat{f}$ or $-\hat{f}$ is a monodromy operator.

We can suppose $\hat f$ is a monodromy operator: if not the following argument holds replacing $\hat f,f$ with $-\hat f,-f$.

By \cite[Corollary 6.4]{monNik}, $\hat{f}$ induces a monodromy operator $\tilde{f}$ on $Y$. We only need to prove that $\tilde{f}$ is $f$ itself. 

To settle this, observe that by definition of $\tilde{f}$, it coincides with $f$ on both $\Lambda_{fix}$ and $\Sigma_Y$. But $\Lambda_{fix}\oplus\Z\Sigma_Y$ has finite index in $\LambdaY$. Therefore, $\tilde{f}=f$.

By the same argument, if $-\hat f$ was a monodromy operator, we would get that there is a monodromy operator $h\in \op{Mon}^2(Y)\subset O^+(\LambdaY)$ such that $h_{|\Lambda_{fix}(2)}=-f_{|\Lambda_{fix}(2)}$  and $h(\Sigma_Y)=f(\Sigma_Y)$. \endproof
\end{proof}

\section*{Acknowledgments}
I would like to thank my supervisors, Christian Lehn and Giovanni Mongardi, for their many useful discussions and suggestions. I am also grateful to Stevell Muller for the many valuable and precise remarks he shared during our email exchange on this problem, in particular for pointing out a missing detail in \cref{extendIsometries}. Finally, I thank Annalisa Grossi for directing me to the relevant references on lattice theory and for her help in navigating them.

\printbibliography
\end{document}